\documentclass{amsart}

\newtheorem{theorem}{Theorem}
\newcommand{\bt}{\begin{theorem}}
\newcommand{\et}{\end{theorem}}
\newtheorem{lemma}{Lemma}
\newcommand{\bl}{\begin{lemma}}
\newcommand{\el}{\end{lemma}}
\newtheorem{corollary}{Corollary}
\newcommand{\bc}{\begin{corollary}}
\newcommand{\ec}{\end{corollary}}
\newtheorem{problem}{Problem}
\newcommand{\bprob}{\begin{problem}}
\newcommand{\eprob}{\end{problem}}
\newtheorem*{conjectureNN}{Conjecture}
\newcommand{\bconjNN}{\begin{conjectureNN}}
\newcommand{\econjNN}{\end{conjectureNN}}
\newtheorem*{theoremNN}{Theorem}
\newcommand{\btNN}{\begin{theoremNN}}
\newcommand{\etNN}{\end{theoremNN}}
\newcommand{\beq}{\begin{equation}}
\newcommand{\eeq}{\end{equation}}
\newcommand{\benum}{\begin{enumerate}}
\newcommand{\eenum}{\end{enumerate}}
\newcommand{\N}{\ensuremath{ \mathbf N }}

\newcommand{\mcf}{\ensuremath{ \mathcal F}}

\newcommand{\mcx}{\ensuremath{ \mathcal X}}

\newcommand{\bmat}{\left(\begin{matrix}}
\newcommand{\emat}{\end{matrix}\right)}

\subjclass[2010]{13E15, 13B02, 13E05, 13F20.} 
\keywords{Polynomial ring, Noetherian ring, intermediate ring, finitely generated ring}

\begin{document}

\title{Finitely generated and  not finitely generated rings}
\author{Melvyn B.  Nathanson}\address{Department of Mathematics\\Lehman College (CUNY)\\Bronx, NY 10468}\email{melvyn.nathanson@lehman.cuny.edu}

\maketitle

\begin{abstract}  Let $R_1$ be a commutative ring, let $R_2$ be a finitely generated extension ring of $R_1$, and let $S$ be a ring that is intermediate between $R_1$ and $R_2$.  For $R_1 = R[x]$ and $R_2 = R[x,y]$, there are simple combinatorial constructions of intermediate rings that are not finitely generated over $R[x]$.  \end{abstract}

Let $R_1$ and $R_2$ be commutative rings with $R_1 \subseteq R_2$.  
The ring $R_2$ is \emph{finitely generated as a ring} over $R_1$ 
if there is a finite subset $\mcx$ of $R_2$ 
such that every element of $R_2$ can be represented 
as a linear combination of monomials in \mcx\ with coefficients in  $R_1$.
The ring $R_2$ is \emph{finitely generated as a module} over $R_1$ 
if there is a finite subset $\mcx$ of $R_2$ 
such that every element of $R_2$ can be represented 
as a linear combination of elements of \mcx\ with coefficients in  $R_1$.

Let  $R_2$ be finitely generated as a ring over $R_1$.  
By Hilbert's basis theorem, if  $R_1$ is Noetherian, then $R_2$ is also Noetherian. 
Let $S$ be a ring that is \emph{intermediate} between $R_1$ and $R_2$, that is,  
\[
R_1 \subseteq S \subseteq R_2.
\]
Artin and Tate~\cite{arti-tate51}  proved that if $R_1$ is Noetherian 
and if $R_2$ is finitely generated as a module over $S$, 
then $S$ is finitely generated as a ring over $R_1$.   
They used this to prove Hilbert's Nullstellensatz
(cf. Zariski~\cite{zari47}, Kunz~\cite[Lemma 3.3]{kunz13}).   

It is natural to ask if \emph{every} intermediate ring $S$ is finitely generated 
as a ring over $R_1$.  
The answer is ``no," and the purpose  of this note is to give simple 
combinatorial constructions  
of intermediate rings $S$ that are not finitely generated over $R_1$.

Let \N\ denote the set of positive integers and $\N_0$ the set of nonnegative integers.  

\bt                                \label{FGR:theorem:main}
Let $\lambda$ be a positive real number or $\lambda = \infty$.  
Let  $\Lambda$ be a subset of $\N \times \N_0$ 
with $(1,0) \in \Lambda$ such that 
\beq                 \label{FGR:sup1}
\sup \left( \frac{b}{a} : (a,b) \in \Lambda  \right) = \lambda 
\eeq
and 
\beq                 \label{FGR:sup2}
\frac{b}{a} < \lambda \qquad \text{for all $ (a,b) \in \Lambda  $.}
\eeq
Consider the set of monomials 
\[
M(\Lambda) =  \left\{  x^ay^b : (a,b) \in \Lambda \right\}.  
\]
Let $R$ be a commutative ring, and let 
$R[M(\Lambda) ]$ be the subring of $R[x,y]$ generated by 
$M(\Lambda)$.  
Then
\[
R[x] \subseteq R[M(\Lambda)] \subseteq R[x,y]
\]
and $R[M(\Lambda)]$ is not finitely generated as a ring over $R[x]$.  
\et

For example, the set $\Lambda_1 = \{ (1,n): n \in \N_0 \}$ satisfies 
conditions~\eqref{FGR:sup1} and~\eqref{FGR:sup2} with $\lambda = \infty$, 
the corresponding set of monomials is 
\[
M(\Lambda_1) = \left\{ x,xy,xy^2, xy^3,\ldots  \right\},
\]
and the ring 
\[
R[M(\Lambda_2)] = R[x,xy,xy^2, xy^3,\ldots ] 
\]
is intermediate between $R[x]$ and $R[x,y]$.   
Similarly, if $(f_n)_{n=-1}^{\infty}$ is the sequence of Fibonacci numbers with 
 $f_{-1} = 1$, $f_0 = 0$, and  $f_1 = 1$, then the set 
 $\Lambda_2 = \{ (f_{2n-1},f_{2n}): n \in \N_0 \}$ satisfies 
conditions~\eqref{FGR:sup1} and~\eqref{FGR:sup2} with $\lambda = (\sqrt{5}+1)/2$.
By Theorem~\ref{FGR:theorem:main}, the intermediate rings $R[M(\Lambda_1)]$ 
and $R[M(\Lambda_2)]$ are not finitely generated over $R[x]$.

Note that inequalities~\eqref{FGR:sup1} and~\eqref{FGR:sup2}  imply that the sets $\Lambda$ 
and $M(\Lambda)$ are infinite.  
If $\lambda = \infty$, then~\eqref{FGR:sup1} implies~\eqref{FGR:sup2}

\begin{proof}
Because $(1,0) \in \Lambda$, we  have $x \in M(\Lambda)$ 
and $R[x] \subseteq R[M(\Lambda)]  \subseteq R[x,y]$.  

Let \mcf\ be a finite subset of $R[M(\Lambda)]$.  
For every polynomial $f$ in $\mcf$, there is a finite 
set $M^*(f)$ of monomials in $M(\Lambda)$ such that $f$  is a linear combination 
of products of monomials in $M^*(f)$.  
This set of monomials is not necessarily unique 
(for example, $(xy)(xy^4) = (xy^2)(xy^3)$), 
but we choose, for each polynomial $f$ in $\mcf$, 
one set $M^*(f)$ of monomials 
in $M(\Lambda)$ that generate $f$.  
Because \mcf\ is a finite set of polynomials, the set  
\[
M^*(\mcf) = \bigcup_{f \in \mcf} M^*(f)  
\]
is a finite set of monomials in $M(\Lambda)$.   
Moreover, $f \in R[M^*(\mcf)]$ for all $f \in \mcf$, and so 
\[
R[\mcf] \subseteq R[M^*(\mcf)] \subseteq R[M(\Lambda)].  
\]
We shall prove that $R[M^*(\mcf)] \neq R[M(\Lambda)]$.    
  
Let 
\[
\beta = \max \left(  \frac{b}{a} :  x^ay^b \in M^*(\mcf) \right).   
\]
Applying inequality~\eqref{FGR:sup2} to the finite set $M^*(\mcf)$, 
we obtain $\beta < \lambda$.   
If $(A,B) \in \N \times \N_0$ and 
$x^A y^B \in R[M^*(\mcf)]$, then $x^A y^B$ is an $R$-linear combination 
of products of monomials in $M^*(\mcf)$.  This implies that 
$x^A y^B$ is a product of monomials in $M^*(\mcf)$.  
Thus, there is a finite sequence $( (a_i,b_i))_{i=1}^n$ 
of ordered pairs in $ \N \times \N_0$ such that 
$x^{a_i}y^{b_i} \in M^*(\mcf)$ for all $i=1,\ldots, n$ and 
\[
x^A y^B = \prod_{i=1}^n x^{a_i}y^{b_i} 
= x^{\sum_{i=1}^n a_i} y^{\sum_{i=1}^n b_i}.
\]  
It follows that  
\[
\frac{B}{A} = \frac{ \sum_{i=1}^n b_i }{ \sum_{i=1}^n a_i} 
\leq\frac{ \beta \sum_{i=1}^n a_i }{ \sum_{i=1}^n a_i}= \beta.
\]
Condition~\eqref{FGR:sup1} implies that the ring $R[M(\Lambda)]$ 
contains monomials $x^A y^B$ with $\beta < B/A < \lambda$.  
It follows that $x^Ay^B \notin R[M^*(\mcf)] $ and so $R[M^*(\mcf)] \neq R[M(\Lambda)]$.   
Therefore,  $R[\mcf] \neq R[M(\Lambda)]$, and 
the ring $R[M(\Lambda)]$ is not finitely generated.    
This completes the proof.  
\end{proof} 

Theorem~\ref{FGR:theorem:main} suggests the following open problems.
Classify the sets $M$ of monomials of the form $x^ay^b$ such that 
\[
R[x] \subseteq R[M] \subseteq R[x,y]
\]
and the ring $R[M]$ is not finitely generated over $R[x]$.  
More generally, describe all  rings $S$ that are intermediate between 
$R[x]$ and $R[x,y]$ and  are not finitely generated over $R$.

I thank Ryan Alweiss for very helpful discussions on this topic at CANT 2016.  

\def\cprime{$'$} \def\cprime{$'$} \def\cprime{$'$}
\providecommand{\bysame}{\leavevmode\hbox to3em{\hrulefill}\thinspace}
\providecommand{\MR}{\relax\ifhmode\unskip\space\fi MR }
\providecommand{\MRhref}[2]{%
  \href{http://www.ams.org/mathscinet-getitem?mr=#1}{#2}
}
\providecommand{\href}[2]{#2}

\end{document}